\documentclass[12pt]{article}
\usepackage{amsmath}
\usepackage{amssymb}
\usepackage[cp1251]{inputenc}
\usepackage[russian]{babel}
\begin{document}

MSC 34C10

\vskip 20pt

\centerline{\bf Oscillatory criteria for the second order linear ordinary differential}

\centerline{\bf   equations  in the marginal sub extremal and extremal cases}

\vskip 10 pt

\centerline{\bf G. A. Grigorian}

\vskip 10 pt

\centerline{0019 Armenia c. Yerevan, str. M. Bagramian 24/5}
\centerline{Institute of Mathematics NAS of Armenia}
\centerline{E - mail: mathphys2@instmath.sci.am}

\vskip 20 pt

\noindent
Abstract. The Riccati equation method is used to establish three new oscillatory criteria for the second order linear ordinary differential equations in the marginal, sub extremal and extremal cases. We show that the first of these criteria  implies the J. Deng's oscillatory criterion. An extremal oscillatory condition for the Mathieu's equation is obtained. The obtained results are compared with some known oscillatory criteria.

\vskip 20pt


\vskip 20pt

1. {\bf Introduction}. Let  $q(t)$  be a real valued locally integrable  function on  $[t_0;+\infty)$. Consider the equation
$$
\phi''(t) + q(t)\phi(t) = 0, \phantom{aaa} t\ge t_0. \eqno (1)
$$

{\bf Definition 1}. The equation (1) is said to be oscillatory, if its
every  solution has arbitrary large zeroes.

Study of the  oscillatory behavior of Eq. (1) is an important problem of the qualitative theory of differential equations and many works are devoted to it (see.  [1] and cited works therein, [2 - 9]). Study of the oscillatory property of Eq. (1) by properties of its coefficients was developed  (and now is being developed) in general in two directions. The goal of the firs direction is to study  the oscillatory property of Eq. (1) on the finite interval (interval oscillatory criteria: Sturm, J. S. W. Wong [2], J. G. Sun, C. H. Ou and J. S. W. Wong [3], Q. Kong [4]).   Then oscillation of Eq. (1) follows from its oscillation on the countable set of finite intervals.  The second one studies correlation between oscillatory property of Eq. (1) and properties of the function $q(t)$ on the whole half axes  (integral oscillatory criteria: Fite, Wintner  [5], Ph. Hartman [3, Theorem 52], Leighton [1], I. V. Kamenev  [6], J. Yan  [7], J. Deng [8], A. Elbert  [9]). Recently  new integral oscillatory criteria were obtained, describing a wide  classes of oscillatory second order linear ordinary differential equations in terms of their coefficients
(see. [5], [8], [9]). Among them note the following important result due to J. Deng   (see [8]) and Q. Kong (see [4]).

{\bf Theorem 1 [8, Theorem 1]}. {\it If for  large  $t\in R$
$$
\int\limits_t^{+\infty}q(\tau)d\tau \ge \frac{\alpha_0}{t} \phantom{a}\mbox{where}\phantom{a} \alpha_0 > \frac{1}{4},
$$
then Eq. (1) is oscillatory.
$\Box$}

{\bf Theorem 2 [4, Theorem 2.3].} {\it Eq. (1) is oscillatory provided that for each $ r \ge t_0$ and for some $\lambda > 1$, either

\noindent
(i) The following two inequalities hold:
$$
\limsup\limits_{t\to+\infty}\frac{1}{t^{\lambda -1}}\int\limits_r^t(t -s)^\lambda q(s) d s > \frac{\lambda^2}{4(\lambda - 1)}
$$
and
$$
\limsup\limits_{t\to+\infty}\frac{1}{t^{\lambda -1}}\int\limits_r^t(s - r)^\lambda q(s) d s > \frac{\lambda^2}{4(\lambda - 1)};\phantom{aaa} \mbox{or}
$$

\noindent
(ii) The following inequality holds
$$
\limsup\limits_{t\to+\infty}\frac{1}{t^{\lambda -1}}\int\limits_r^t(t -s)^\lambda [q(s) + q(2t - s)] d s > \frac{\lambda^2}{2(\lambda - 1)}.
$$
$\Box$}

In this work we use the Riccati equation method for establishing three new oscillatory criteria for Eq. (1)
in three different directions having relations to the cases when: a) $\lim\limits_{t\to \infty}\frac{1}{t}\int\limits_{t_0}^t  d\tau \int\limits_{t_0}^\tau q(s) d s = \lambda \ne \pm\infty$ (marginal case, see below Theorem 3 and Corollary~ 1); b) $\liminf\limits_{t\to \infty}\frac{1}{t}\int\limits_{t_0}^t  d\tau \int\limits_{t_0}^\tau q(s) d s = - \infty$ (sub extremal case, see below Theorem~ 4 and Corollary 2); c) $\lim\limits_{t\to \infty}\frac{1}{t}\int\limits_{t_0}^t  d\tau \int\limits_{t_0}^\tau q(s) d s = - \infty$ ( extremal case, see  below Theorem 5 and Corollary 3). Note that the oscillation of Eq. (1) in the remaining (main) cases when:\linebreak d) $\lim\limits_{t\to \infty}\frac{1}{t}\int\limits_{t_0}^t  d\tau \int\limits_{t_0}^\tau q(s) d s = + \infty$ (regular  case); e) $- \infty < \liminf\limits_{t\to \infty}\frac{1}{t}\int\limits_{t_0}^t  d\tau \int\limits_{t_0}^\tau q(s) d s < \linebreak \limsup\limits_{t\to \infty}\frac{1}{t}\int\limits_{t_0}^t  d\tau \int\limits_{t_0}^\tau q(s) d s$ (irregular  case) was discovered by A. Wintner  (see [2, Theorem~ 51]) and Ph. Hartman (see [2, Theorem~ 52]) respectively.
We show that Corollary 1 (see below) implies Theorem 1. We  obtain an oscillatory condition for the Mathieu's equation. This condition contains extremal cases for the Matheou's equation. We compare the obtained results with  some known oscillatory criteria.

\vskip 20pt

{\bf 2. Main results}.
 Denote by  $\Omega$ the set of all positive and absolutely  continuous  on  $[t_0;+\infty)$  functions $f(t)$ such that the functions $f'(t)^2$ are locally integrable  on  $[t_0;+\infty)$ .

{\bf Theorem 3}. {\it Let for some  $f\in \Omega, \phantom{a} \lambda\in R, \phantom{a} \alpha \ge 1$   the following conditions be satisfied:

\noindent
1)
$
\int\limits_{t_0}^{+\infty}\exp\biggl\{ \int\limits_{t_0}^t \frac{d\tau}{f(\tau)}\int\limits_{t_0}^\tau\biggl[2 f(s) q(s) - \frac{f'(s)^2}{2 f(s)}\biggr] d s\biggr\} d t
=+\infty
$;

\noindent
2) $\liminf \limits_{t\to +\infty}\bigl\{ \frac{1}{f(t)}\int\limits_{t_0}^t \bigl[ 4 f(\tau) q(\tau) - \frac{f'(\tau)^2}{f(\tau)}\bigr] d \tau - 4 \int\limits_{t_0}^t q(\tau) d \tau\bigr\} < +\infty$;

\noindent
3) $\int\limits_{t_0}^{+\infty}\exp\bigl\{ - 4 \lambda t + 4 \int\limits_{t_0}^t d\tau \int\limits_{t_0}^\tau q(s) d s\bigr\} d t < +\infty;$

\noindent
4) $\limsup\limits_{t \to +\infty} \frac{1}{t^\alpha}\int\limits_{t_0}^t(t - \tau)^\alpha q(\tau) d \tau \ge \lambda.$.

\noindent
Then Eq. (1) is oscillatory. $\Box$}

{\bf Remark 1}. For  $f(t) \equiv 1$ the condition 2) of Theorem 3 is  always satisfied.

{\bf Remark 2}. The condition 4) of Theorem 3 can be replaced by one of the following conditions
$$
\int\limits_{t_0}^{+\infty}\biggl[\lambda - \int\limits_{t_0}^\tau q(s) d s \biggr]^2 d \tau < +\infty, \phantom{aaa} \int\limits_{A_\lambda} d \tau = + \infty,
$$
where  $A_\lambda \equiv \bigl\{ t\ge t_0 : \int\limits_{t_0}^t q(\tau) d \tau  \ge \lambda\bigr\}$.
Indeed, when one of these conditions takes place we have: $\int\limits_{t_0}^{+\infty}\bigl[y_*(\tau) + \lambda - \int\limits_{t_0}^\tau q(s) d s \bigr]^2 d \tau = +\infty.$ From here and from (14) it follows (13) (see below the proof of Theorem 3).

{\bf Corollary  1}. {\it Let for some $\lambda \in R$  the following conditions be satisfied:

\noindent
A) $\lim\limits_{t \to +\infty} \frac{1}{t}\int\limits_{t_0}^t d \tau \int\limits_{t_0}^\tau q(s) d s = \lambda$;

\noindent
B)
$\frac{1}{t}\int\limits_{t_0}^t d \tau \int\limits_{t_0}^\tau q(s) d s  \le \lambda - \frac{\ln Q_0(t)}{4 t} + \frac{Q_1(t)}{t},$

\noindent
for enough large   $t$,  where  $Q_k(t)\phantom{a} ( k =0,1)$ are some real continuous functions on  $[t_0;+\infty)$, satisfying the conditions:
$Q_0(t) > 0, \phantom{a}  t\ge t_0; \phantom{a} \int\limits_{t_0}^{+\infty}
\frac{d\tau}{Q_0(\tau)} < +\infty, \phantom{a} \sup\limits_{t\ge t_0} Q_1(t) < +\infty.$

\noindent
Then Eq. (1) is oscillatory. $\Box$}

Let us show that  Theorem 1 is a consequence of Corollary~1. Suppose the condition of Theorem 1 holds. Then from the convergence of the integral  $\int\limits_{t_0}^{+\infty} q(\tau) d \tau$ it follows existence of the finite limit
$$
\lim\limits_{t\to +\infty} \frac{1}{t}\int\limits_{t_0}^t d\tau \int\limits_{t_0}^\tau q(s) d s  \stackrel{def} = \lambda. \eqno (2)
$$
(Since by L'Hospital's rule $\lim\limits_{t\to +\infty} \frac{1}{t}\int\limits_{t_0}^t
d\tau \int\limits_\tau^{+\infty} q(s) d s = \lim\limits_{t\to +\infty} \int\limits_t^{+\infty} q(s) d s = 0$). Therefore \linebreak condition  $A)$ of Corollary 1 takes place. According to the condition of Theorem 1 chose  $\eta_0 > 0$  so large that   $\int\limits_t^{+\infty} q(\tau) d \tau \ge \frac{\alpha_0}{t}, \phantom{a}  t\ge \eta_0$.  Then taking into account (2) we will have:
$\frac{1}{t}\int\limits_{t_0}^t d\tau \int\limits_{t_0}^\tau q(s) d s  = \frac{1}{t}\int\limits_{t_0}^{\eta_0} d\tau \int\limits_{t_0}^\tau q(s) d s + \frac{1}{t}\int\limits_{\eta_0}^t d\tau \int\limits_{t_0}^{+\infty} q(s) d s - \frac{1}{t}\int\limits_{\eta_0}^t d\tau \int\limits_\tau^{+\infty} q(s) d s \le  \lambda  - \linebreak - \frac {\alpha_0 \ln t }{t} + \frac{c_1}{t}, \phantom{a} t\ge \eta_0,$
where  $c_1 \equiv \int\limits_{t_0}^{\eta_0} d\tau \int\limits_{t_0}^\tau q(s) d s - \lambda \eta_0 + \alpha_0 \ln \eta_0$.  Therefore condition $B)$  of Corollary 1 takes place. Thus the Theorem 1 is a consequence of Corollary 1. The next example shows that Corollary~1 does not follow (even with additional hypothesis of convergence of integral  $\int\limits_{t_0}^{+\infty} q(\tau) d\tau$)  from Theorem 1.

{\bf Example 1}. Consider the equation
$$
\phi''(t) + \biggl[\frac{\alpha_0}{t^2} + \frac{\alpha \cos(\beta t)}{t^\gamma}\biggr]\phi(t) = 0, \phantom{aaa} t\ge t_0 > 0. \eqno (3)
$$
where   $\alpha_0, \phantom{a} \alpha, \phantom{a} \beta, \phantom{a} \gamma $ are some real constants, $\alpha_0 > \frac{1}{4}, \phantom{a} \gamma > 0, \phantom{a} \alpha \beta \ne 0$.   Without loss of generality we will assume that

$\sin(\beta t_0) = \frac{\gamma}{\beta} \frac{\cos(\beta t_0)}{t_0} - \gamma(\gamma + 1) t_0^\gamma \int\limits_{t_0}^\infty\frac{\cos(\beta \tau)}{\tau^{\gamma +2}} d \tau\stackrel{def}{=} \xi(t_0)$ (since $\xi(t_0) \to 0$ for $t_0 \to \infty$).  Then
$$
\int\limits_{t_0}^t \biggl[\frac{\alpha_0}{\tau^2} + \frac{\alpha\cos(\beta\tau)}{\tau^\gamma}\biggr] d\tau = \frac{\alpha_0}{t_0} - \frac{\alpha_0}{t}   + \frac{\alpha \sin(\beta t)}{\beta t^\gamma} -
 \frac{\alpha}{\beta}\frac{\sin(\beta t_0)}{t_0^\gamma} + \frac{\alpha\gamma}{\beta}\int\limits_{t_0}^t\frac{\sin(\beta \tau)}{\tau^{\gamma+1}}  d \tau =
 $$
$$
=\frac{\alpha_0}{t_0} - \frac{\alpha_0}{t}   + \frac{\alpha \sin(\beta t)}{\beta t^\gamma} -
 \frac{\alpha}{\beta}\frac{\sin(\beta t_0)}{t_0^\gamma} + \frac{\alpha\gamma}{\beta}\int\limits_{t_0}^\infty\frac{\sin(\beta \tau)}{\tau^{\gamma+1}}  d \tau -
 \frac{\alpha\gamma}{\beta} \int\limits_t^{\infty}\frac{\sin(\beta \tau)}{\tau^{\gamma+1}} d \tau =
$$
$$
=\frac{\alpha_0}{t_0} - \frac{\alpha_0}{t}   + \frac{\alpha \sin(\beta t)}{\beta t^\gamma} -
 \frac{\alpha\gamma}{\beta} \int\limits_t^{\infty}\frac{\sin(\beta \tau)}{\tau^{\gamma+1}} d \tau  - \phantom{aaaaaaaaaaaaaaaaaaaaaaaaaaaaaaaa}
$$
$$
 - \frac{\alpha}{\beta}\frac{\sin(\beta t_0)}{t_0^\gamma} + \frac{\alpha \gamma}{\beta^2}\frac{\cos(\beta t_0)}{t_0^{\gamma+1}}
 - \frac{\alpha\gamma(\gamma+1)}{\beta^2}\int\limits_{t_0}^\infty \frac{\cos(\beta \tau)}{\tau^{\gamma +2}}d\tau =
$$
$$
\phantom{aaaaaaaaaaaaaaaaaaaaaaaaaaaaaaaa}=\frac{\alpha_0}{t_0} - \frac{\alpha_0}{t}   + \frac{\alpha \sin(\beta t)}{\beta t^\gamma} -
 \frac{\alpha\gamma}{\beta} \int\limits_t^{\infty}\frac{\sin(\beta \tau)}{\tau^{\gamma+1}} d \tau, \phantom{aaa} t\ge t_0.
$$
Therefore
$$
\int\limits_{t_0}^t d\tau \int\limits_{t_0}^\tau\biggl[ \frac{\alpha_0}{s^2} + \frac{\alpha \cos(\beta s)}{s^\gamma}\biggr] d s = \frac{\alpha_0 (t - t_0)}{t_0} - \alpha_0 \ln t + \alpha_0 \ln t_0 + \phantom{aaaaaaaaaaaaaaaaaaaaaaaaaaa}
$$
$$
\phantom{aaaaaaaaaaaaaaaaaaaaaaaaaaa}+\frac{\alpha}{\beta} \int\limits_{t_0}^t \frac{\sin(\beta \tau)}{\tau^\gamma} d\tau - \frac{\alpha \gamma}{\beta}\int\limits_{t_0}^t d\tau \int\limits_\tau ^{\infty} \frac{\sin(\beta s)}{s^{\gamma+1}} d s. \eqno (4)
$$
Since $\frac{\sin(\beta t)}{t^\gamma} - \gamma \int\limits_t^{\infty}\frac{\sin(\beta \tau)}{\tau^{\gamma+1}} d \tau \to 0$ for $t \to \infty$ by L'Hospital's rule
$$
\lim\limits_{t\to \infty}\biggl[\frac{1}{t}\int\limits_{t_0}^t\frac{\sin(\beta \tau)}{\tau^\gamma}d \tau - \frac{\gamma}{t}\int\limits_{t_0}^t d\tau \int\limits_\tau^{\infty}\frac{\sin(\beta s)}{s^{\gamma+1}} d s\biggr] = 0
$$
From here and from (4) it follows  $\lim\limits_{t\to \infty} \frac{1}{t} \int\limits_{t_0}^t d\tau\int\limits_{t_0}^\tau \bigl[\frac{\alpha_0}{s^2} + \frac{\alpha \cos(\beta s)}{s^\gamma}\bigr] d s = \frac{\alpha_0}{t_0} \stackrel{def} = \lambda$.
 Therefore condition  $A)$ of Corollary 1 is satisfied. From (4) we have
$$
\frac{1}{t}\int\limits_{t_0}^t d\tau\int\limits_{t_0}^\tau\biggl[ \frac{\alpha_0}{s^2} + \frac{\alpha \cos(\beta s)}{s^\gamma}\biggr] d s = \lambda - \alpha_0 \frac{\ln t}{t} + \alpha_0 \frac{\ln t_0}{t} + \phantom{aaaaaaaaaaaaaaaaaaaaaaaaaaaaaaa}
$$
$$
\phantom{aaaaaaaaaaaaaaaaaaaaaaaaaaa}+\frac{\alpha}{\beta t} \int\limits_{t_0}^t \frac{\sin(\beta \tau)}{\tau^\gamma} d\tau - \frac{\alpha \gamma}{\beta t}\int\limits_{t_0}^t d\tau \int\limits_\tau ^{+\infty} \frac{\sin(\beta s)}{s^{\gamma+1}} d s, \phantom{aaa} t\ge t_0.
$$
From here it follows that for Eq.  (3) condition  $B)$ of Corollary 1 is satisfied too. So Eq. (3) is oscillatory. It is easy to show that
$$
\int\limits_t^{+\infty}\biggl[\frac{\alpha_0}{\tau^2} + \frac{\alpha\cos(\beta\tau)}{\tau^\gamma}\biggr] d\tau = \frac{\alpha_0}{t} -  \frac{\alpha \sin ( \beta t)}{ \beta t^\gamma} + \frac{O(t)}{t^{\gamma + 1}}, \phantom{aaa} t\ge t_0,
$$
where $O(t)$ is a bounded function on $[t_0;+\infty)$. Hence it is clear that for  $0< \gamma < 1$  as well as for  $ \gamma =1, \phantom{a} \bigl|\frac{\alpha}{\beta}\bigr| > \alpha_0$ Theorem 1 is not applicable to Eq. (3).

{\bf Remark 3}.   {\it The conditions of Corollary 1 exclude the condition of Ph. Hartman's oscillatory criterion  [2, Theorem 52] and the condition of I. V. Kamenev's  oscillatory criterion.}

{\bf Example 2}. Consider the equation
$$
\phi''(t) + \biggl[\frac{1}{4t^2} + \frac{1}{4t^2 ln t} + \dots + \frac{1 + \varepsilon}{4 t^2 ln t \dots ln_r t} + \frac{\alpha\sin(\beta t)}{t ln t}\biggr] \phi(t) = 0, \phantom{aaa} ln_r t > 0. \eqno (5)
$$
Here $ln_1t = ln t, \dots ln_r t = ln ln_{r-1} t, \phantom{a} r=2, 3, \dots, \phantom{a} \varepsilon, \phantom{a} \alpha, \phantom{a} \beta$ are some real constants, $\varepsilon > 0, \phantom{a} \beta \ne 0$. Using the identity
$\frac{1}{t}\int\limits_{t_0}^t d \tau \int\limits_{t_0}^\tau q(s) d s = \int\limits_{t_0}^t q(s) d s - \frac{1}{t}\int\limits_{t_0}^t s q(s) d s, \phantom{a} t\ge t_0$, one can readily check that for Eq. (5) all conditions of Corollary 1 are fulfilled. Therefore Eq. (5) is oscillatory. It is not difficult to verify that Theorem 1 and Theorem 2 are not applicable to Eq. (5).

{\bf Theorem  4}. {\it  Let for some  $\varepsilon > 0,\phantom{a}\alpha \ge 1$  and  $\lambda\in R $                         the following conditions be satisfied:

\noindent
5) $\int\limits_{t_0}^{+\infty}\exp\biggr\{ - 4\lambda t + 4\int\limits_{t_0}^td\tau\int\limits_{t_0}^\tau q(s) d s\biggr\} d t = +\infty;$

\noindent
6) $\liminf\limits_{t \to +\infty} \frac{1}{t^\alpha}\int\limits_{t_0}^t(t - \tau)^\alpha q(\tau) d \tau \le \lambda - \varepsilon$.

\noindent
Then Eq. (1) is oscillatory. $\Box$}

Denote: $B_\lambda \equiv \{t \ge t_0: \int\limits_{t_0}^t d \tau \int\limits_{t_0}^\tau q(s) d s \ge \lambda t\}$.

{\bf Corollary 2}. {\it  Let for some  $\lambda \in R$   and   $\varepsilon > 0$ the following conditions be satisfied:

\noindent
C) $\int\limits_{B_\lambda} d \tau = +\infty$;

\noindent
D) $\liminf \limits_{t\to + \infty} \frac{1}{t}\int\limits_{t_0}^t d\tau \int\limits_{t_0}^\tau q(s) d s \le \lambda  - \varepsilon$.

\noindent
Then Eq. (1) is oscillatory. $\Box$}

{\bf Remark 4}. From the conditions of Corollary 2  is seen that
$$
\liminf\limits_{t\to +\infty}\frac{1}{t}\int\limits_{t_0}^t d\tau \int\limits_{t_0}^\tau q(s) d s < \limsup\limits_{t\to +\infty}\frac{1}{t}\int\limits_{t_0}^t d\tau \int\limits_{t_0}^\tau q(s) d s
$$
However Corollary 2 is not a consequence of Ph. Hartman's oscillatory criterion \linebreak [2. Theorem  52]. Indeed in particular case when
$$
\int\limits_{t_0}^t d\tau \int\limits_{t_0}^\tau q(s) d s = \left\{
\begin{array}{l}
t\sin^3 t, \phantom{a} \sin t \ge 0;\\
t^2\sin^3 t, \phantom{a} \sin t < 0,
\end{array}\right.
$$
$ t\ge t_0$, for  $\lambda = 0, \phantom{a} \varepsilon =1$
all conditions of Corollary  2 are satisfied, whereas for this case the conditions of Ph. Hartman's criterion are not fulfilled  (since  $\liminf\limits_{t\to +\infty}\frac{1}{t}\int\limits_{t_0}^t d\tau \int\limits_{t_0}^\tau q(s) d s =~-\infty$).

{\bf Theorem 5}. {\it Let for some $f \in \Omega$ the condition 1) and the condition

\noindent
7) \phantom{a} $\liminf\limits_{ t \to +\infty}\biggl\{\frac{1}{f(t)}\int\limits_{t_0}^t\bigl[2 f(\tau) q(\tau) - \frac{f'(\tau)^2}{2 f(\tau)}\bigr]d \tau - 2\int\limits_{t_0}^t q(\tau) d \tau \biggr\} = - \infty$

\noindent
be satisfied. Then Eq. (1) is oscillatory.$\Box$}

Consider the Matheu's equation (see [10], p. 111)
$$
\phi''(t) + (\delta + \varepsilon \cos (2t)) \phi(t) = 0, \phantom{aaa} t\ge t_0. \eqno (6)
$$
Here $\delta$ and $\varepsilon$ are some real constants, $\varepsilon \ne 0$. Set: $F_\varepsilon (\mu) \equiv - \frac{\pi |\varepsilon| \mu}{4(\pi + \mu)} + \frac{\mu^2}{\pi + \mu} \int\limits_{-\frac{\pi}{4}}^{\frac{\pi}{4}}\frac{\sin^2(2t)}{1 + \mu\cos (2t)} d t, \linebreak \mu \ge 0, \phantom{aa} m(\varepsilon) \equiv \inf\limits_{\mu \ge 0} F_\varepsilon(\mu).$ Obviously $m(\varepsilon) < 0$.

{\bf Corollary 3}. {\it If $\delta > m(\varepsilon)$    then Eq. (6) is oscillatory. $\Box$}

{\bf Remark 5}. {\it For $\delta \in (m(\varepsilon); 0)$ we have the extremal case of Eq. (1) : $\int\limits_{t_0}^{+\infty} q(t) d t = -\infty$. In this case Theorem 2 is not applicable to Eq. (6).}

{\bf Remark 6}. {\it Using standard methods for integration of trigonometric functions it is easy to calculate the value of the integral, presenting  in the expression  for $F_\varepsilon(\mu)$:
$$
\frac{\mu^2}{\pi + \mu} \int\limits_{-\frac{\pi}{4}}^{\frac{\pi}{4}}\frac{\sin^2(2t)}{1 + \mu\cos (2t)} d t  = \left\{
\begin{array}{l}
\frac{\pi- 2\mu}{2(\pi + \mu)} - \frac{2\sqrt{1 - \mu^2}}{\pi + \mu} \arctan\sqrt{\frac{1 - \mu}{1 + \mu}}, \phantom{aaa} 0 \le \mu <1;\\
\phantom{a}\\
\frac{\pi - 2}{2(\pi + 1)}, \phantom{aaaaaaaaaaaaaaaaaaaaaaaa} \mu =1;\\
\phantom{a}\\
\frac{\pi- 2\mu}{2(\pi + \mu)}  + \frac{\sqrt{\mu^2 - 1}}{\pi + \mu} \ln \frac{\sqrt{\mu + 1} + \sqrt{\mu - 1}} {\sqrt{\mu + 1} + \sqrt{\mu - 1}}, \phantom{aaa} \mu > 1.
\end{array}
\right.
$$
}

{\bf Example 3}. It is not difficult to verify that $m(4) < F_4(1) = - \frac{\pi + 2}{2(\pi +1)}$. Therefore by Corollary 3 the equation
$$
\phi''(t) + \biggl[ - \frac{\pi + 2}{2(\pi +1)} + 4\cos (2t)\biggr] \phi(t) = 0, \phantom{aaa} t\ge t_0. \eqno (7)
$$
is oscillatory.

{\bf Remark 7}. It is not difficult to verify that the Sturm's comparison criterion (compari- \linebreak son of Eq. (7) with an equation $\phi''(t) + q_0\phi(t) = 0$, where $q_0$ is any constant $> 0$) is not applicable to Eq. (7).

One can readily  verify that the criteria of Ph. Hartman [2, Theorem 52] and I. V. Kamenev  [6] are not applicable to Eq. (3). For  $\gamma \ge 1$  the J. Yan's criterion  [7, p. 277, THEOREM] is not applicable  to Eq.  (3) too.  None of Theorems 1 - 9 of work  [5] is applicable to Eq.  (3), and for  $\gamma \ge \frac{1}{2}$  Theorem 11 of work [5] does not applicable  to Eq. (3) as well. Theorem 1 as well as the criteria of Ph. Hartman [2, Theorem 52], I. V. Kamenev [6], and J. Yan [7, p. THEOREM] are not applicable to Eq. (7).

{\bf 3. Proof of the main results}. To prove the main results we need in some auxiliary propositions.

{\bf 3.1. Auxiliary propositions}.  Let  $a(t)$  and $b(t)$  be real valued locally integrable functions on  $[t_0;+\infty)$. Consider the Riccati equation
$$
y'(t) + y^2(t) + a(t) y(t) + b(t) = 0, \phantom{aaa} t\ge t_0. \eqno (8)
$$

{\bf Definition 2}. A solution of Eq. (8) is said to be  $t_1$-regular, if it exists on  $[t_1;+\infty)$ $(t_1 \ge t_0)$.

{\bf Definition  3}.A $t_1$-regular solution  $y(t)$  of Eq. (8) is said to be  $t_1$-normal, if there exists $\delta > 0$ such that every solution $y_1(t)$ of Eq. (8) with  $y_1(t_1) \in (y(t_1) - \delta; y(t_1) + \delta)$     is $t_1$-regular, otherwise it is said to be  $t_1$-extremal.

Denote by  $reg(t_1)$  the set of all  $y_{(0)} \in R$,
for which the solution  $y(t)$ of Eq. (8) with  $y_(t_1) = y_{(0)}$  is $t_1$-regular.

{\bf Lemma 1}. {\it If Eq. (8) has a $t_1$-regular solution then it has the unique  $t_1$-extremal solution  $y_*(t)$,  moreover  $reg(t_1) = [y_*(t_1); +\infty)$.}

This lemma is proved in  [11] for the case of continuous  $a(t)$  and $b(t)$. For the case of locally integrable  $a(t)$  and $b(t)$ the proof by analogy.

{\bf Lemma 2}. {\it Suppose Eq.  (8) has a $t_1$-regular solution and $b(t) \ge 0, \phantom{a} t\ge t_0$.  Then if:

$I) \int\limits_{t_0}^{+\infty}\exp\biggl\{ - 2\int\limits_{t_0}^t a(\tau) d\tau\biggr\} d t< +\infty$, then $ y_*(t) < 0, \phantom{a} t\ge t_2$  for some  $t_2 \ge t_1$;

$ II)\int\limits_{t_0}^{+\infty}\exp\biggl\{ - 2\int\limits_{t_0}^t a(\tau) d\tau\biggr\} d t= +\infty$, then
$y_*(t) \ge 0, \phantom{a} t\ge t_1$,

\noindent
where  $y_*(t)$ is the unique  $t_1$-extremal solution of Eq.  (8).}

This lemma is proved in  [12] for the case of continuous  $a(t)$  and $b(t)$. For the case of locally integrable  $a(t)$  and $b(t)$ the proof by analogy.

Let  $y(t)$  be a  $t_1$-regular solution of Eq. (8). Consider the integral
$$
\nu_y(t) \equiv \int\limits_t^{+\infty}\exp\biggl\{ - \int\limits_t^\tau \bigl[2 y(s) + a(s)\bigr]d s\biggr\}d t, \phantom{aaa} t\ge t_1.
$$

{\bf Theorem 6 [11, Theorem 2.A]}.{\it The integral  $\nu_y(t)$ converges for all  $t\ge t_1$ if and only if  $y(t)$  is   $t_1$-normal. $\Box$}

This theorem is proved in  [11] for the case of continuous  $a(t)$  and $b(t)$. For the case of locally integrable  $a(t)$  and $b(t)$ the proof by analogy.

{\bf 3.2. Proof of rhe main results.    Proof of Theorem 3}.
  Suppose Eq. (1)  is not oscillatory. Then the equation
$$
x'(t) + x^2(t) + q(t) = 0, \phantom{aaa} t \ge t_0, \eqno (9)
$$
has a  $t_1$-regular solution for some  $t_1 \ge t_0$ (see [13], p. 332).
In this equation make the substitution
$$
x(t) = y(t) + \lambda - \int\limits_{t_0}^t q(\tau)d\tau, \phantom{aaa} t\ge t_0. \eqno (10)
$$
We obtain:
$$
y'(t) + y^2(t) + 2\biggl(\lambda - \int\limits_{t_0}^t q(\tau) d\tau\biggr) y(t) + \biggl(\lambda - \int\limits_{t_0}^t q(\tau) d\tau\biggr)^2 = 0,\phantom{aaa} t\ge t_0. \eqno (11)
$$
Since Eq. (9) has a  $t_1$-regular solution, from  (10) it follows that the last equation has a $t_1$-regular solution too.
Then since   $\biggl(\lambda - \int\limits_{t_0}^t q(\tau) d\tau\biggr)^2 \ge 0, \phantom{a} t\ge t_0,$  by virtue of Lemma  2.I) from condition 3) it follows that
$$
y_*(t) < 0, \phantom{aaa}t \ge t_2, \eqno (12)
$$
for some $t_2 \ge t_1$, where $y_*(t)$ is the unique  $t_1$-extremal solution of Eq. (11).
Let us show that
$$
y_*(t) \to -\infty \phantom{a}\mbox{for} \phantom{a} t \to +\infty \eqno (13)
$$
By virtue of  (11) we have:
$$
y_*(t) = y_*(t_1) - \int\limits_{t_1}^t\biggl[y_*(\tau) + \lambda - \int\limits_{t_0}^\tau q(s) d s\biggr]^2 d\tau, \phantom{aaa} t\ge t_1. \eqno (14)
$$
Suppose that the relation (13) is not true. Then from (12) and (14) it follows that $y^*(t)$ is a decreasing function on  $[t_1;+\infty)$ with a negative finite limit:
$$
y_*(+\infty)\equiv \lim\limits_{t\to +\infty} y_*(t) < 0 \phantom{aaa}(y_*(t)\downarrow y_*(+\infty) > - \infty) \eqno (15)
$$
From here and (14) it follows that  $\int\limits_{t_1}^{+\infty}\biggl[t_*(\tau) + \lambda - \int\limits_{t_0}^\tau q(s) d s\biggr]^2 d\tau < +\infty$. Then
$$
0 \le \lim \limits_{t\to +\infty} \frac{1}{t^\alpha}\int\limits_{t_1}^t\biggl(t - \tau\biggr)^{\alpha - 1}\biggl[y_*(\tau) + \lambda - \int\limits_{t_0}^\tau q(s) d s\biggr]^2 d\tau \le \phantom{aaaaaaaaaaaaaaaaaaaaaaaaaaaaaaaaaaa}
$$
$$
\phantom{aaaaaaaaaaaaaaaaaa}\le \lim \limits_{t\to +\infty} \frac{1}{t^\alpha}\int\limits_{t_1}^{+\infty} \biggl(t - \tau\biggr)^{\alpha - 1}\biggl[y_*(\tau) + \lambda - \int\limits_{t_0}^\tau q(s) d s\biggr]^2 d\tau = 0. \eqno (16)
$$
Set: $\rho(t) \equiv y_*(t) - y_*(+\infty), \phantom{a} t\ge t_1;\phantom{aa}I\equiv \limsup\limits_{t \to +\infty}\frac{1}{t^\alpha}\int\limits_{t_1}^t \bigl( t - \tau\bigr)^{\alpha -1}\bigl[y_*(\tau) + \lambda  - \int\limits_{t_0}^\tau q(s) d s\bigr]^2 d \tau.$

\noindent
It is evident that
$$
\rho(t) \to 0 \phantom{a} for \phantom{a} t \to +\infty. \eqno (17)
$$
We have:
$$
I = \limsup\limits_{t \to +\infty}\frac{1}{t^\alpha}\int\limits_{t_1}^t\bigl(t - \tau)^{\alpha - 1}\biggl[y_*(+\infty) + \lambda - \int\limits_{t_0}^\tau q(s) d s + \rho(\tau)\biggr]^2 d\tau = \phantom{aaaaaaaaaaaaaaaaaaaaaa}
$$
$$
= \limsup\limits_{t \to +\infty}\biggl[\frac{y_*^2(+\infty)}{t^\alpha}\int\limits_{t_1}^t \bigl(t - \tau\bigr)^{\alpha - 1} d \tau + \frac{2 y_*(+\infty)}{t^\alpha}\int\limits_{t_1}^t \bigl(t - \tau\bigr)^{\alpha - 1} \rho(\tau) d \tau +
$$
$$
 +\frac{2 y_*(+\infty)}{t^\alpha} \int\limits_{t_1}^t \bigl(t - \tau\bigr)^{\alpha - 1}\biggl(\lambda - \int\limits_{t_0}^\tau q(s) d s\biggr) d \tau + \frac{1}{t^\alpha}\int\limits_{t_1}^t \biggl(\lambda - \int\limits_{t_0}^\tau q(s) d s + \rho(\tau)\biggr)^2 d\tau\biggr]\ge
$$
$$
\ge \limsup\limits_{t \to +\infty}\biggl[\frac{y_*^2(+\infty)}{t^\alpha}\int\limits_{t_1}^t \bigl(t - \tau\bigr)^{\alpha - 1} d \tau + \frac{2 y_*(+\infty)}{t^\alpha}\int\limits_{t_1}^t \bigl(t - \tau\bigr)^{\alpha - 1} \rho(\tau) d \tau +
$$
$$
 +\frac{2 y_*(+\infty)}{t^\alpha} \int\limits_{t_1}^t \bigl(t - \tau\bigr)^{\alpha - 1}\biggl(\lambda - \int\limits_{t_0}^\tau q(s) d s\biggr) d \tau\biggr]. \eqno (18)
$$
From the condition 4) and from  (15) it follows that
$$
\lim\limits_{n \to +\infty} \frac{2y_*(+\infty)}{\theta_n^\alpha}\int\limits_{t_1}^{\theta_n}\bigl(\theta_n - \tau\bigr)^{\alpha - 1}\biggl(\lambda - \int \limits_{t_0}^\tau q(s) d s\biggr)d \tau \ge 0, \eqno (19)
$$
for some infinitely large sequence  $\{\theta_n\}_{n=1}^{+\infty}$.
Obviously  by virtue of (17) \linebreak $\lim\limits_{t\to +\infty} \frac{1}{t^\alpha}\int\limits_{t_1}^t \bigl(t - \tau\bigr)^{\alpha - 1} \rho(\tau) d \tau = 0.$ From here, from (18) and (19) it follows that \linebreak $I \ge \frac{y_*^2(+\infty)}{\alpha}>~0,$ which contradicts  (16).
The obtained contradiction proves (13).  It follows from the condition 2) that there exists a infinitely large sequence  $\{\xi_n|\}_{n=1}^{+\infty}$ such that $ S\equiv \sup\limits_{n\ge 1}\biggl\{\frac{1}{f(\xi_n)}\int\limits_{t_0}^{\xi_n}\bigl[4f(\tau)q(\tau) - \frac{f'(\tau)^2}{f(\tau)}\bigl] d\tau - 4 \int\limits_{t_0}^{\xi_n} q(\tau) d\tau\biggr\} < +\infty$.
In view of this and relation (9) we chose  $t_3 = \xi_{n_0}$  so large that
$$
y_*(t_3) + \lambda + S/4 < 0. \eqno (20)
$$
Show that the solution  $x_0(t)$ of Eq. (9) with
$$
x_0(t_3) = \frac{1}{f(t_3)}\int\limits_{t_0}^{t_3}\biggl[\frac{f'(\tau)^2}{4 f(\tau)} - f(\tau) q(\tau)\biggr] d \tau   \eqno (21)
$$
is $t_3$-normal.  By  5) from (20) it follows  $x_*(t_3) = y_*(t_3) - \int\limits_{t_0}^{t_3} q(\tau) d \tau \le x_0(t_3)$. By virtue of Lemma 1 from here it follows that  $x_0(t)$ is $t_3$-normal. Since $x_0(t)$ is a  $t_3$-regular solution of Eq. (9), we have
$$
f(t) x_0'(t) +  f(t)x_0^2(t) + f(t) q(t) = 0, \phantom{aaa} t\ge t_3.
$$
Let us integrate this equality from  $t_3$  to $t$. We  get
$$
f(t)x_0(t) = \int\limits_{t_3}^t\bigl[f(\tau) x_0^2(\tau) - f'(\tau)x_0(\tau)\bigr]d\tau = f(t_3) x_0(t_3) - \int\limits_{t_3}^t f(\tau) q(\tau) d \tau, \phantom{aaa} t \ge t_3.
$$
After completing the square under the integral on the left-hand side of this equality and dividing both sides of it on  $f(t)$
we obtain
$$
x_0(t) = \frac{1}{f(t)}\int\limits_{t_3}^t f(\tau) \biggl[x_0(\tau) - \frac{f'(\tau)}{2f(\tau)}\biggr]^2 d \tau = \frac{c}{f(t)}  - \frac{1}{f(t)}\int\limits_{t_0}^t \biggl[\frac{f'(\tau)^2}{f(\tau)} - f(\tau) q(\tau)\biggr] d \tau, \eqno (22)
$$
where  $c\equiv f(t_3) x_0(t_3) - \int\limits_{t_0}^{t_3}\bigl[\frac{f'(\tau)^2}{4 f(\tau)} - f(\tau) q(\tau)\bigr]d \tau$. By virtue of (21) we have $c = 0$. Therefore from (22) we get: $-2 x(t) \ge \frac{1}{f(t)}\int\limits_{t_0}^t \bigl[ 2 f(\tau) q(\tau) - \frac{f'(\tau)^2}{f(\tau)}\bigr] d \tau, \phantom{a} t\ge t_3$. Then
$$
\int\limits_{t_3}^{+\infty}\exp\biggl\{ -2 \int\limits_{t_3}^t x_0(\tau) d \tau\biggr\} d t \ge M\int\limits_{t_3}^{+\infty}\exp\biggl\{\int\limits_{t_3}^t \frac{d\tau}{f(\tau)}\int\limits_{t_0}^\tau\biggl[2 f(s) q(s) - \frac{f'(s)^2}{2 f(s)}\biggr] d s\biggr\} d t, \eqno (23)
$$
where   $M \equiv \exp\bigl\{ - \int\limits_{t_0}^{t_3}\frac{d\tau}{f(\tau)}\int\limits_{t_0}^\tau\bigl[2 f(s) q(s) - \frac{f'(s)^2}{2 f(s)}\bigr] d s\bigr\} = const > 0.$ Since  $x_0(t)$ is  $t_3$-normal,  by virtue of Theorem 2  the left-hand side of the inequality  (23) is finite. Whereas from the condition 1) it follows, that its right-hand side equals to  $+\infty$. The obtained contradiction proves the theorem.

{\bf Proof of Corollary 1}. At first we prove the corollary in the particular case when $\lambda > 0$.   Take: $f(t) \equiv 1$.  Then from  А) it follow 1) and 4), and by virtue of Remark 1 the condition  2) is fulfilled. From В) it follows:
$$
\int\limits_{t_0}^{+\infty}\exp\biggl\{ - 4\lambda + 4\int\limits_{t_0}^t d\tau \int\limits_{t_0}^\tau q(s) d s \biggr\} d t \le   c_0 +  \int\limits_{\eta_0}^{+\infty} \exp\biggl\{ - t\biggl| 4\lambda    - \frac{4}{t}\int\limits_{t_0}^t d\tau \int\limits_{t_0}^\tau q(s) d s\biggr|\biggr\} d t \le
$$
$$
\le c_0  + \int\limits_{\eta_0}^{+\infty} \exp\biggl\{ - \ln Q_0(t) + Q_1(t)\biggr\} d t \le c_0 + M \int\limits_{\eta_0}^{+\infty} \frac{d\tau}{Q_0(\tau)} < +\infty,
$$
where
$
c_0 \equiv  \int\limits_{t_0}^{\eta_0} \exp\biggl\{ - t\biggl| 4\lambda    - \frac{4}{t}\int\limits_{t_0}^t
d\tau
\int\limits_{t_0}^\tau q(s) d s\biggr|\biggr\} d t, \phantom{a} M \equiv \exp\{\sup\limits_{t\ge t_0}
Q_1(t)\} < +\infty,$ and
 \phantom{a}$ \eta_0
$  is a  enough large number such that for all  $t \ge \eta_0$  the inequality of condition  В) holds. So for $\lambda > 0$   all conditions of Theorem 3 are fulfilled, and for this particular case the corollary is proved.
Suppose  $\lambda \le 0$. Consider the equation
$$
\phi''(t) + \widetilde{q}(t) \phi(t) = 0, \phantom{aaa} t  \ge t_0. \eqno (24)
$$
where  $\widetilde{q}(t) = q(t) + \Delta q(t), \phantom{a} \Delta q(t)$ is a continuous function on  $[t_0; +\infty)$  such that  $\Delta q(t) \ge~0, \linebreak  t\in [t_0; t_0 + 1], \phantom{a} \Delta q(t) = 0, \phantom{a} t\ge t_0 + 1, \phantom{a} \int\limits_{t_0}^{t_0 + 1} \Delta q(\tau) d\tau  > |\lambda|.$                            Then
$$
\lim\limits_{t\to +\infty}\frac{1}{t}\int\limits_{t_0}^t d \tau \int\limits_{t_0}^\tau \widetilde{q}(s) d s =
\lim\limits_{t\to +\infty}\biggl[\frac{1}{t}\int\limits_{t_0}^t d \tau \int\limits_{t_0}^\tau q(s) d s + \frac{1}{t}\int\limits_{t_0}^{t_0 +1} d \tau \int\limits_{t_0}^\tau \Delta q(s) d s + \phantom{aaaaaaaaaaaa}
$$
$$
\phantom{aaaaaaaaaaaaaaaaaaaaaaaaaaaaaaaaaa} + \frac{1}{t}\int\limits_{t_0+1}^t d \tau \int\limits_{t_0}^\tau \Delta q(s) d s\biggr] = \lambda + \int\limits_{t_0}^{t_0+1} \Delta q(s) d s > 0.
$$
Thus

\noindent
$\widetilde{A} ) \phantom{a} \lim\limits_{t\to +\infty}\frac{1}{t}\int\limits_{t_0}^t d \tau \int\limits_{t_0}^\tau \widetilde{q}(s) d s = \widetilde{\lambda} > 0$

\noindent
From  В) it follows:
$$
\lim\limits_{t\to +\infty} \frac{1}{t}\int\limits_{t_0}^t d\tau \int\limits_{t_0}^\tau\widetilde{q}(s) d s = \frac{1}{t}\int\limits_{t_0}^t d\tau \int\limits_{t_0}^\tau q(s) d s  + \frac{1}{t}\int\limits_{t_0}^t d\tau \int\limits_{t_0}^\tau \Delta q(s) d s \le
$$
$$
\le \lambda  - \frac{\ln Q_0(t)}{4 t} + \frac{Q_1(t)}{t} +  \frac{1}{t}\int\limits_{t_0}^t d\tau \int\limits_{t_0}^\tau \Delta q(s) d s
$$
for all enough large  $t$.   From here it follows:

\noindent
$
\widetilde{B}) \phantom{a} \frac{1}{t}\int\limits_{t_0}^t d\tau \int\limits_{t_0}^\tau\widetilde{q}(s) d s \le \widetilde{\lambda}  - \frac{\ln Q_0(t)}{4t} + \frac{\widetilde{Q}_1(t)}{t}
$

\noindent
where  $\widetilde{Q}_1(t) \equiv Q_1(t) + \int\limits_{t_0}^{t_0 + 1} d \tau \int\limits_{t_0}^\tau \Delta q(s) d s - (t_0 + 1) \int\limits_{t_0}^{t_0 + 1}\Delta q(s) d s$.
Obviously  $\sup\limits_{t\ge t_0} \widetilde{Q}_1(t) < + \infty$. Then by already proven from $\widetilde{A})$    and $\widetilde{B})$ it follows  that Eq. (20) is oscillatory.
Since   $\widetilde{Q}(t)$   differs from  $q(t)$  only at most on  $[t_0; t_0 + 1]$, from the oscillation  of Eq. (20) it follows the oscillation of Eq. (1). The corollary is proved.

{\bf Proof of Theorem 4}.  Let  $y_*(t)$  be the same  $t_1$-extremal solution of Eq. (11), as in the proof of Theorem~3. By virtue of Lemma 2. II) from nonnegativity of the function \linebreak $\Bigl(\lambda - \int\limits_{t_0}^t q(\tau) d \tau\Bigr)^2, \phantom{a} t\ge t_0,$  and from  5)
it follows that  $y_*(t)$ is a nonnegative decreasing function  with finite limit:
$$
y_*(+\infty) \equiv \lim\limits_{t\to +\infty} y_*(t) \ge 0 \phantom{aaa} \bigl( y_*(t) \downarrow y_*(+\infty) \ge 0\bigr). \eqno (25)
$$
From here and from (14) it follows, that   $\int\limits_{t_1}^{+\infty}\bigl[y_*(\tau) + \lambda - \int\limits_{t_0}^\tau q(s) d s\bigr]^2 d\tau < +\infty.$  Then
$$
\lim\limits_{t\to +\infty} \frac{1}{t^\alpha}\int\limits_{t_1}^t \bigl(t - \tau\bigr)^{\alpha - 1}\biggl[y_*(\tau) + \lambda - \int\limits_{t_0}^\tau q(s) d s\biggr]^2 d\tau = 0. \eqno (26)
$$
Set (as above):  $\rho(t) \equiv y_*(t) - y_*(+\infty), \phantom{a} t \ge t_1.$  Then we have
$$
\limsup\limits_{t\to +\infty}\frac{1}{t^\alpha}\int\limits_{t_1}^t\bigl(t - \tau\bigr)^{\alpha -1}\biggl[y_*(\tau) + \lambda - \int\limits_{t_0}^\tau q(s) d s\biggr]^2 d\tau = \phantom{aaaaaaaaaaaaaaaaaaaaaaaaaaaaaaa}
$$
$$
= \limsup\limits_{t\to +\infty}\biggl[\frac{(y_*(+\infty) + \varepsilon)^2}{t^\alpha}\int\limits_{t_1}^t\bigl(t - \tau\bigr)^{\alpha -1} d\tau + \frac{2 (y_*(+\infty) + \varepsilon)}{t^\alpha}\int\limits_{t_1}^t\bigl(t - \tau\bigr)^{\alpha -1}\bigl(\lambda -\phantom{aaaa}
 $$
$$
\phantom{aaaaaaaaaa} -\varepsilon - \int\limits_{t_0}^\tau q(s) d s\bigr) d \tau + \frac{1}{t^\alpha}\int\limits_{t_1}^t\bigl(t - \tau\bigr)^{\alpha -1}\bigl(\lambda - \varepsilon - \int\limits_{t_0}^\tau q(s) d s + \rho (\tau)\bigr)^2 d \tau\biggr] \ge
$$
$$
\ge \limsup\limits_{t\to +\infty} \biggl[\frac{(y_*(+\infty) + \varepsilon)^2}{t^\alpha}   \int\limits_{t_1}^t\bigl(t - \tau\bigr)^{\alpha -1}d\tau + \frac{2 (y_*(+\infty) + \varepsilon)}{t^\alpha}\int\limits_{t_1}^t\bigl(t - \tau\bigr)^{\alpha -1}\rho(\tau) d \tau +\phantom{aaaaaaa}
$$
$$
\phantom{aaaaaaaaaaaaaaaa}+\frac{2 (y_*(+\infty) + \varepsilon)}{t^\alpha}\int\limits_{t_1}^t\bigl(t - \tau\bigr)^{\alpha -1}\bigl(\lambda - \varepsilon - \int\limits_{t_0}^\tau q(s) d s\bigr)d\tau\biggr]. \eqno (27)
$$
Since  $\rho(t) \to 0$   for  $t\to +\infty$,  we have
$$
\lim\limits_{t\to +\infty} \frac{1}{t^\alpha}\int\limits_{t_1}^t\bigl(t - \tau\bigr)^{\alpha -1}\rho(\tau) d \tau = 0. \eqno (28)
$$
From the condition 6) it follows, that there exists an infinitely large sequence  $\{t_n\}_{n=2}^{+\infty}$ such that $$
\lim\limits_{n \to +\infty}\frac{1}{t_n^\alpha}\int\limits_{t_1}^{t_n}\bigl(t_n - \tau\bigr)^{\alpha - 1} d \tau \int\limits_{t_0}^\tau q(s) d s \le \frac{\lambda - \varepsilon}{\alpha}.
$$
From here, from  (27) and  (28) it follows that
$$
\limsup\limits_{t\to +\infty}\frac{1}{t^\alpha}\int\limits_{t_1}^t\bigl(t - \tau\bigr)^{\alpha -1}\biggl[y_*(\tau) + \lambda - \int\limits_{t_0}^\tau q(s) d s\biggr]^2 d\tau  \ge  \frac{[y_*(+\infty) + \varepsilon]^2}{\alpha} > 0,
$$
which contradicts   (26). The obtained contradiction proves the theorem.

{\bf Proof of Corollary 2}. From the condition  С) it follows
$$
\int\limits_{t_0}^{+\infty}\exp\biggl\{ - 4\lambda t + 4 \int\limits_{t_0}^t d\tau\int\limits_{t_0}^\tau q(s) d s\biggr\} d t \ge \int\limits_{B_\lambda}exp\biggl\{ - 4\lambda t + 4 \int\limits_{t_0}^t d\tau\int\limits_{t_0}^\tau q(s) d s\biggr\} d t = +\infty.
$$
Therefore the condition 5) of Theorem 4 is satisfied. It follows from D) that for  $\alpha =1$  the condition  6) of Theorem 4 is satisfied too. The corollary is proved.

{\bf Proof of Theorem 5}.  Suppose Eq. (1) is not oscillatory. Then Eq. (9)
has a $t_1$-regular solution for some $t_1 \ge t_0$. Then by (10) the equation
$$
y'(t) + y^2(t) + 2\Bigl(\lambda - \int\limits_{t_0}^t q(\tau) d \tau\Bigr) y(t) + \Bigl(\lambda - \int\limits_{t_0}^t q(\tau) d \tau\Bigr)^2 = 0, \phantom{aaa} t\ge t_0. \eqno (29)
$$
has also a $t_1$ - regular solution. Let then $y_*(t)$ be the $t_1$ - extremal solution of Eq. (29). We have
$$
y_*(t) = y_*(t_1) - \int\limits_{t_1}^t\Bigl[y_*(\tau) + \lambda - \int\limits_{t_0}^\tau q(s) d s\Bigr]^2 d\tau, \phantom{aaa} t \ge t_1.
$$
Due to the equality 7) chose $t_3 \ge t_1$ such that
$$
\frac{1}{f(t_3)}\int\limits_{t_0}^{t_3}\Bigl[f(\tau) q(\tau) - \frac{f'(\tau)^2}{4 f(\tau)}\Bigr] d\tau - \int\limits_{t_0}^{t_3} q(\tau) d \tau < - \lambda - y_*(t_1). \eqno (30)
$$
Let $x_0(t)$ be the solution of Eq. (28) with
$$
x_0(t_3) = \frac{1}{f(t_3)}\int\limits_{t_0}^{t_3}\biggr[\frac{f'(\tau)^2}{4f(\tau)} - f(\tau) q(\tau)\biggr]d\tau.
$$
Show that $x_0(t)$ is $t_3$-normal. By (30) we have:
$$
x_*(t_3) = y_*(t_3) = \lambda - \int\limits_{t_0}^{t_3} q(\tau) d\tau \le  y_*(t_1) = \lambda - \int\limits_{t_0}^{t_3} q(\tau) d\tau <
$$
$$
< - \frac{1}{f(t_3)}\int\limits_{t_0}^{t_3}\biggr[f(\tau) q(\tau) -  \frac{f'(\tau)^2}{4f(\tau)}\biggr]d\tau = x_0(t_3).
$$
Therefore $x_0(t)$ is $t_3$ normal. Further the continuation of the proof  is similar to the proof of Theorem 3. The theorem is proved.

{\bf Proof of Corollary 3}. It is not difficult to verify that $F_\varepsilon(\mu)$ reaches its minimum on $(0;+\infty)$. Let then $m(\varepsilon) = F_\varepsilon (\mu_0)$ for some $\mu_0 > 0$ and let $g(t)$ be a periodic function of  period $\pi$, defined by formula
$$
g(t)=\left\{
\begin{array}{l}
1 + \mu_0 \cos 2t, \phantom{a} t\in [- \frac{\pi}{4};\frac{\pi}{4}];\\
\phantom{a}\\
1, \phantom{aaaaaaaaaaa}t\in (\frac{\pi}{4};\frac{3 \pi}{4}].
\end{array}
\right.
$$
From the condition of Corollary 3 it follows that
$$
\int\limits_\frac{-\pi}{4}^\frac{\pi}{4}\biggl[2 g(t)(\delta + \varepsilon \cos 2t) - \frac{g'(t)^2}{g(t)}\biggr] d t > 0, \eqno (31)
$$
Consider the sequence of functions: $h_1(t), \dots, h_k(t), \dots,$ defined by formulae
$$
h_k(t)\equiv\left\{
\begin{array}{l}
k(t+ \frac{1}{k}), \phantom{aaaa} k \in [-\frac{1}{k};0];\\
\phantom{a}\\
-k(t- \frac{1}{k}), \phantom{aaa} k \in (0;\frac{1}{k}];\\
\phantom{a}\\
0, \phantom{aaaaaaaaaa} k \not\in [-\frac{1}{k};\frac{1}{k}],
\end{array}
\right.
$$
$k=1,2,\dots$. It is easy to show that
$$
\int\limits_{-\frac{1}{k}}^\frac{1}{k}\frac{h_k'(t)^2}{1 + h_k(t)} d t > k, \phantom{aaa} k=1,2,\dots. \eqno (32)
$$
Let $k_0 < k_1 < \dots $ be  a  increasing sequence of  natural numbers (yet arbitrary). Define the function $f(t)$ on $[-\frac{\pi}{4}; + \infty)$ as follows
$$
f(t) = g(t) + h_{k_{j+1}}\Bigl(t + \frac{\pi}{4} - \pi k_j + \frac{2}{k_{j+1}}\Bigr), \phantom{aa} t\in \biggl[- \frac{\pi}{4} + \pi k_{j-1}; - \frac{-\pi}{4} + \pi k_j\biggr], \phantom{aa} j=1,2,\dots.
$$
Using (31) we chose $k_1$ so large that
$$
\int\limits_{-\frac{\pi}{4}}^t \frac{d\tau}{f(\tau)}\int\limits_{-\frac{\pi}{4}}^\tau\biggl[2f(s)(\delta + \varepsilon \cos 2 s) - \frac{f'(s)^2}{2 f(s)}\biggr] d s > 0, \phantom{a} t \in \Bigl[- \frac{\pi}{4} + \pi(k_1 -1); - \frac{\pi}{4} + \pi k_1- \frac{2}{k_1}\Bigr],
$$
and using (32) chose $k_2 > k_1$ so large that
$$
\frac{1}{f(-\frac{\pi}{4} + \pi k_1)}\int\limits_{-\frac{\pi}{4}}^{-\frac{\pi}{4} + \pi k_1}\biggl[2 f(\tau)(\delta + \varepsilon \cos 2\tau) - \frac{f'(\tau)^2}{2 f(\tau)}\biggr] d \tau  - \int\limits_{-\frac{\pi}{4}}^{-\frac{\pi}{4} + \pi k_1} (\delta + \varepsilon \cos 2\tau)d \tau  < -1.
$$

Using (31) we chose $k_3 > k_2$ so large that
$$
\int\limits_{-\frac{\pi}{4}}^t \frac{d\tau}{f(\tau)}\int\limits_{-\frac{\pi}{4}}^\tau\biggl[2f(s)(\delta + \varepsilon \cos 2 s) - \frac{f'(s)^2}{2 f(s)}\biggr] d s > 0, \phantom{a} t \in \Bigl[- \frac{\pi}{4} + \pi(k_3 -1); - \frac{\pi}{4} + \pi k_3- \frac{2}{k_3}\Bigr],
$$
and using (32) chose $k_4 > k_3$ so large that
$$
\frac{1}{f(-\frac{\pi}{4} + \pi k_3)}\int\limits_{-\frac{\pi}{4}}^{-\frac{\pi}{4} + \pi k_3}\biggl[2 f(\tau)(\delta + \varepsilon \cos 2\tau) - \frac{f'(\tau)^2}{2 f(\tau)}\biggr] d \tau  - \int\limits_{-\frac{\pi}{4}}^{-\frac{\pi}{4} + \pi k_3} (\delta + \varepsilon \cos 2\tau)d \tau  < -2,
$$
so on .... We see that for such a chose of $k_1 < k_2 < \dots$ for the constructed function $f(t)$ all the conditions of Theorem 5 for Eq. (6) are fulfilled. The corollary is proved

\vskip 20 pt

\centerline{\bf References}

\vskip 12pt
\noindent
1. C. A.  Swanson,  Comparison and oscillation theory of linear differential equations.\linebreak  \phantom{a}  Academic press, New York and London, 1968.

\noindent
2. Q. Kong, M. Pasic, Second Order Differential Equations: Some Significant Results Due\\\phantom{a} to James S. W. Wong. Differential Equations and Applications, Vol. 6. Number 1(2014) \\\phantom{a}99 - 163.

\noindent
3. J. G. Sun, C. H. Ou and J. S. W. Wong, Interval Oscillation Theorems for a Second\linebreak  \phantom{aa} Order Linear Differential Equations. Computers and Mathematics with   Applicatins.\linebreak  \phantom{aa}  48 (2004) 1693 - 1699.

\noindent
4. Q. Kong, Interval Criteria for Oscillation of Second-Order Linear Ordinary Differential \linebreak \phantom{a} Equations, Journal of Mathematical analysis and applications, 229, 258 - 270 (1999).

\noindent
5. M. K. Kwong, Integral criteria for second- order linear oscillation. Electronic journal \linebreak  \phantom{a} of   qualitative theory of differential equations, 2006, $N^\circ$ 10, pp. 1 - 28.\linebreak \phantom{a} http// WWW. math.u- szeged.hu/ejqtde/

\noindent
6. I. V. Kamenev,  On oscillation criterion for the second order linear differential equations, \\\phantom{a}  Matematicheskie zametki, vol 23, $N^\circ$ 2, 1978, pp. 249 - 251.

\noindent
7. J. Yan, Oscillation theorems for second order linear differential equations with damping. \\\phantom{a} Proceedings of the American Mathematical Society, vol. 98, $N^\circ$  2, 1986,  pp. 276 - 282.

\noindent
8. J. Deng, Oscillation criteria for second order linear differential equations. Journal of  \\\phantom{aa}Mathematical Analysis  and Applications 271(2002) 283 - 287.

\noindent
9. A. Elbert, Oscillation/Nonoscillation Criteria for Linear Second Order Differential\linebreak  \phantom{a} Equations. Journal ournal of Mathematical Analysis  and Applications, 226,\linebreak  \phantom{a} 207 - 219 (1998).

\noindent
10. L. Cezari, Asymptotic Behavior and Stability Problems in Ordinary Differential \linebreak  \phantom{aa}  equations, Berlin 1959. Translated under the title Asimptoticheskoe povedenie i \linebreak  \phantom{aa} ustoichivost' reshenii obiknovennich differentialnikh urqvnenii, Moskow, 1964.

\noindent
11. G. A. Grigorian, Properties of Solutions of Riccati Equation. Journal of Contemporary \linebreak  \phantom{aa}Mathematical analysis. 2007, vol. 42, No. 4, pp.184 - 197.

\noindent
12.  G. A. Grigorian, On the Stability of Systems of Two First-Order Linear Ordinary \linebreak  \phantom{aa} Differential Equations. Differential Equations (Original Russian Differential'nie Urav-\linebreak \phantom{aa} nenia) 2015, Vol. 51, No 3, pp. 1 -10.

\noindent
13. Ph. Hartman, Ordinary differential equations, Second edition, The Jhon Hopkins \linebreak  \phantom{aa}University, Baltimore, Merilend, 1982.

\end{document}